\def \NN {\mathbb N}
\def \CC {\mathbb C}
\def \QQ {\mathbb Q}
\def \RR {\mathbb R}
\def \ZZ {\mathbb Z}

\def \epsilon{\varepsilon}

\def \A  {{\mathcal A}}

\def \D  {{\mathcal D}}
\def \H  {{\mathcal H}}

\def \O  {{\mathcal O}}

\def \S  {{\mathcal S}}

\def \d {\text{d}}
\def \FF {\Bbb F}

\def \M {{\mathcal M}}

\def \bfalpha {{\boldsymbol{\alpha}}}

\def \si {\sigma}

\newcommand{\res}{\text{res}}

\newcommand{\lsum}{\sum\limits}
\newcommand{\lprod}{\prod\limits}

\documentclass[12pt,reqno]{amsart}

\usepackage{amsmath,amssymb}

\usepackage{url} 

\begin{document}


\title[]{Converse theorems: from the Riemann zeta function \\ to the Selberg class}
\author[]{ALBERTO PERELLI}
\maketitle

{\bf Abstract.} This is an expanded version of the author's lecture at the XX Congresso U.M.I., held in Siena in September 2015. After a brief review of $L$-functions, we turn to the classical converse theorems of H.Hamburger, E.Hecke and A.Weil, and to some later developments. Finally we present several results on converse theorems in the framework of the Selberg class of $L$-functions.



\vskip.5cm
\section{$L$-functions}

\smallskip
An $L$-{\it function} is, loosely speaking, a Dirichlet series associated with a structure $\A$, encoding a good amount of information on $\A$. Classically, such a structure may be of arithmetic, geometric or automorphic nature; in this section we briefly describe several important examples. A Dirichlet series is a function of the complex variable $s=\si+it$, defined in some right half-plane as
\[
f(s) = \sum_{n=1}^\infty \frac{a_n}{n^s} \hskip2cm a_n\in\CC.
\]
The link with the structure $\A$ is obtained by a suitable choice of the coefficients $a_n$ or, in many cases, by expressing $f(s)$ as an Euler product of type
\[
f(s) = \prod_p f_p(s),
\]
where $p$ runs over the prime numbers and the local factors $f_p(s)$ are defined in  terms of $\A$. Once the $L$-function is set up, the first problem is to study its analytic properties (meromorphic continuation and polar structure, functional equation, distribution of zeros, order of growth,...). Then, one faces the problem of deriving interesting information on $\A$ from such analytic properties. Both steps are highly non-trivial and, even in the simplest cases, require good ideas coupled with solid techniques. This is, in broad terms, Riemann's approach to the study of prime numbers. Since their first appearance on the scene of prime numbers, $L$-functions have seen a constant growth of their range of applications, and have undergone an intense study of their analytic properties; however, many central problems are still open. The state of the art is well expressed by Bombieri's words \cite[p.38]{Bom/2010}: ``{\it Zeta and $L$-functions, starting from the mysterious first object found by Riemann, have by now grown into an even more mysterious set of objects that encode the most delicate aspects and questions of arithmetic.}'' 

\medskip
In the literature one often finds the name {\it zeta functions} to denote some of the objects which here are called by the generic name $L$-functions. This is essentially due to historical reasons; we shall keep the original names when dealing with classical examples of such functions. Moreover, we shall present $L$-functions in {\it normalized form}, namely with the point $s=1/2$ as center of the critical strip.

\bigskip
{\bf 1.1. Riemann zeta function.} We start with Riemann's mysterious first object, namely the {\it Riemann zeta function},  ancestor of all $L$-functions. It is defined for $\sigma>1$ by the absolutely convergent Dirichlet series
\[
\zeta(s) := \sum_{n=1}^\infty \frac{1}{n^s}
\]
and, thanks to the unique factorization in $\ZZ$, it has there the Euler product expansion
\[
\zeta(s) = \prod_p\big(1-\frac{1}{p^s}\big)^{-1}.
\]
Actually, the Euler product of $\zeta(s)$ is an analytic equivalent of unique factorization. The identity
\begin{equation}
\label{1-1}
\sum_{n=1}^\infty \frac{1}{n^s}  = \prod_p(1-\frac{1}{p^s})^{-1},
\end{equation}
called Euler identity and discovered for real values of $s$ by Euler in 1737, is the starting point of analytic number theory. Indeed, the right hand side offers a strong link to the prime numbers, while the left hand side can be studied by the methods of complex analysis. In 1859, Riemann obtained the basic properties of his zeta function, namely: $\zeta(s)$ can be extended by analytic continuation to a meromorphic function over $\CC$, its only singularity is a simple pole at $s=1$ with residue 1 and satisfies the functional equation
\begin{equation}
\label{1-new1}
\pi^{-s/2} \Gamma\big(\frac{s}{2}\big)\zeta(s) = \pi^{-(1-s)/2} \Gamma\big(\frac{1-s}{2}\big)\zeta(1-s),
\end{equation}
where $\Gamma(s)$ is Euler's function. The functional equation provides a kind of symmetry of the values of $\zeta(s)$ with respect to the critical line $\si=1/2$. This applies in particular to the non-trivial zeros of $\zeta(s)$, lying in the critical strip $0\leq \si\leq 1$, which are symmetric with respect to such a line. Riemann's ideas are a masterpiece of insight and elegance, and inspired a huge development of $L$-functions theory. 

\medskip
Euler used \eqref{1-1} to give a new proof of the existence of infinitely many primes, completely different from Euclid's one, by letting $s\to1^+$ and observing that the left hand side tends to $\infty$. Having \eqref{1-1} for complex $s$, Riemann was able to get a much more explicit link between $\pi(x)$, the number of primes up to $x$, and $\zeta(s)$. Such a link is better expressed by means of the Mellin transform, a variant of the Fourier transform, in the form
\begin{equation}
\label{1-2}
\psi(x) = \frac{1}{2\pi i} \int_{2-i\infty}^{2+i\infty} -\frac{\zeta'(s)}{\zeta(s)} \frac{x^s}{s}\d s,
\end{equation}
where $\psi(x)$ is a kind of weighted version of $\pi(x)$ and $-\zeta'(s)/\zeta(s)$ is the logarithmic derivative of $\zeta(s)$ (rigorously, \eqref{1-2} requires that $x$ is not a prime power, and the integral is meant as principal value). Since the transition from $\psi(x)$ to $\pi(x)$ (and viceversa) is easy, \eqref{1-2} represents a formidable starting point for investigating the asymptotic behavior of $\pi(x)$ as $x\to\infty$. Indeed, one may shift the line of integration to the left, crossing the simple poles of $-\zeta'(s)/\zeta(s)$ at $s=1$ and at the zeros of $\zeta(s)$ (and of $1/s$ at $s=0$), thus getting by the residue theorem that
\begin{equation}
\label{1-3}
\psi(x) = x - \sum_\rho \frac{x^\rho}{\rho} + \ \text{small error,}
\end{equation}
where $\rho$ runs over the non-trivial zeros of $\zeta(s)$. Since the sum over $\rho$ is only conditionally convergent, formula \eqref{1-3}, called the {\it explicit formula}, is therefore not really good for a direct treatment, but suitable approximated versions are available. However, it is already clear from \eqref{1-3} that the location of the zeta zeros has a direct impact on the asymptotic behavior of $\psi(x)$, i.e. on the distribution of prime numbers. Indeed, the contribution of a non-trivial zero $\rho=\beta+i\gamma$ to $\psi(x)$ is roughly $x^\beta/|\gamma|$, hence the zeros with larger real part somehow dominate the sum in \eqref{1-3}.

\medskip
The Euler product immediately shows that $\zeta(s)\neq0$ for $\si>1$, and Riemann's ideas led J.Hadamard and Ch.J. de la Vall\'ee-Poussin in 1896 to the proof of the prime number theorem
\[
\pi(x) \sim \frac{x}{\log x},
\]
conjectured by Gauss around 1792 and by Legendre in 1808, when they proved that $\zeta(1+it)\neq 0$ for every $t\in\RR$. Actually, de la Vall\'ee-Poussin obtained the zero-free region $\si>1-c/\log(|t|+2)$, with some $c>0$, which implies a more precise asymptotic formula for $\pi(x)$. The best known zero-free region of $\zeta(s)$, due to I.M.Vinogradov and his school in the 1950's, is only slightly sharper than de la Vall\'ee-Poussin's one, although its proof requires much harder work. Most probably, the zeta zeros are distributed according to the famous {\it Riemann Hypothesis} (RH for short), asserting that all non-trivial zeros lie on the critical line. RH implies the optimal distribution of prime numbers, in the sense that the asymptotic formula for $\pi(x)$ becomes
\begin{equation}
\label{1-4}
\pi(x) = \text{Li}(x) + O(x^{1/2}\log x);
\end{equation}
here Li$(x)$ denotes the logarithmic integral
\[
\text{Li}(x) := \int_2^x\frac{\d t}{\log t} = \frac{x}{\log x} + \frac{x}{\log^2x} + \dots + \frac{(N-2)!x}{\log^{N-1}x} + O\big(\frac{x}{\log^Nx}\big)
\]
for any given integer $N\geq 2$. Actually, \eqref{1-4} is equivalent to RH; moreover, it is known unconditionally that the remainder term in \eqref{1-4} has oscillations of size, roughly, $x^{1/2}$.

\medskip
The Riemann Hypothesis is probably the most famous open problem in mathematics; again with Bombieri's words \cite[p.111]{Bom/2006}: ``{\it The failure of the Riemann hypothesis would create havoc in the distribution of prime numbers. This fact alone singles out the Riemann hypothesis as the main open question of prime number theory.}'' There is rather strong evidence in favor of RH. First of all, numerical evidence: the first $10^{13}$ zeros lie on the critical line and the same
holds true for billions of zeros around heights $10^{14}$, $10^{15}$, \dots, $10^{24}$. Moreover, there are several results inside the Riemann zeta function theory supporting RH, notably the existence of a positive proportion of zeros lying on the critical line, proved by A.Selberg in 1942, while every half-plane $\si>\alpha$ with $\alpha>1/2$ contains only an infinitesimal proportion of them. Finally, there is a good amount of evidence coming from analogies with other theories, such as the Riemann Hypothesis for the zeta functions of varieties over finite fields, and the consistency of certain statistical properties of the zeta zeros and of the eigenvalues of random matrices.

\medskip
We refer the interested reader to the standard textbook of Ingham \cite{Ing/1990} for an introduction to the Riemann zeta function and its applications to prime number theory, to Titchmarsh's treatise \cite{Tit/1986} for a comprehensive presentation of the finer theory of $\zeta(s)$ and to Bombieri's survey \cite{Bom/2010} of the classical theory of zeta functions.

\bigskip
{\bf 1.2. $L$-functions associated with algebraic number fields.} The natural generalization of the Riemann zeta function to an algebraic number field $K$ is the {\it Dedekind zeta function}, defined for $\si>1$ by the absolutely convergent Dirichlet series
\[
\zeta_K(s) := \sum_{\frak{A}\neq0} \frac{1}{N(\frak{A})^s} = \sum_{n=1}^\infty \frac{a(n)}{n^s}.
\]
Here $\frak{A}$ runs over the non-zero ideals of the ring of integers $\O_K$ of $K$, $N(\frak{A})$ is the norm of $\frak{A}$ and $a(n)$ is the number of $\frak{A}$ with $N(\frak{A})=n$; in particular, $\zeta(s)=\zeta_\QQ(s)$. Thanks to the unique factorization into prime ideals $\frak{P}$, $\zeta_K(s)$ has the Euler product
\[
\zeta_K(s) = \prod_{\frak{P}} \big(1-\frac{1}{N(\frak{P})^s}\big)^{-1} = \prod_p\prod_{j=1}^d \big(1-\frac{\alpha_j(p)}{p^s}\big)^{-1},
\]
where $d$ is the degree of $K$, $|\alpha_j(p)|\leq 1$ for all $p$ and $j$, and $|\alpha_j(p)|=1$, $j=1\dots,d$, for all but finitely many $p$. The basic analytic properties of $\zeta_K(s)$ were obtained by E.Hecke in 1917, by a truly remarkable extension of Riemann's techniques to the framework of number fields: $\zeta_K(s)$ can be extended by analytic continuation to a meromorphic function over $\CC$ with a simple pole at $s=1$ with residue $\nu_K$, and satisfies the functional equation
\[
\Lambda_K(s) = \Lambda_K(1-s),
\]
where
\[
\Lambda_K(s) = \big(\frac{\sqrt{|D_K|}}{2^{r_2}\pi^{d/2}}\big)^s\Gamma\big(\frac{s}{2}\big)^{r_1} \Gamma(s)^{r_2} \zeta_K(s), 
\]
$r_1$ and $r_2$ denote the number of real and complex immersions of $K$, respectively, and $D_K$ is the discriminant of $K$. Therefore, $\zeta_K(s)$ has a general analytic structure similar to $\zeta(s)$.

\medskip
Clearly, the properties above, and many others obtainable by suitably extending to $\zeta_K(s)$ the techniques devised for $\zeta(s)$, can be used to derive results on the distribution of prime ideals. Moreover, several classical conjectures about $\zeta(s)$ can be formulated for $\zeta_K(s)$ as well. In other words, the theory of the Dedekind zeta functions, although at present weaker when finer properties are concerned, is in several aspects parallel to that of $\zeta(s)$. However, the greater arithmetical complexity of number fields with respect to the rational field suggests that the functions $\zeta_K(s)$ contain further information on the invariants of $K$. In particular, the value of the residue at $s=1$ is
\[
\nu_K = \frac{h_K2^{r_1}(2\pi)^{r_2}R_K}{\sqrt{|D_K|}w_K},
\]
where $R_K$ and $\omega_K$ are the regulator and the number of roots of unity of $K$, respectively, and $h_K$ is the class number of $K$. The class number, namely the number of ideal classes modulo the principal ideals, is very interesting since $h_K=1$ if and only if $\O_K$ has unique factorization; hence $h_K$ is a kind of measure of the non-unique factorization of the ring $\O_K$. The analytic theory of the functons $\zeta_K(s)$, coupled with standard facts from algebraic number theory, offers a tool for the study of $h_K$, giving fruitful results. However, it is known that $\zeta_K(s)$ does not contain all the information on $K$. Indeed, there exist non-isomorphic number fields with the same Dedekind zeta function; for example, $\zeta_K(s)$ does not determine the class number of $K$. This is a rather common phenomenon: $L$-functions encode much information about the underlying structures $\A$, but in general do not determine $\A$ uniquely.

\medskip
Another important issue about Dedekind zeta functions, which has no counterpart in the theory of the Riemann zeta function, is the factorization of $\zeta_K(s)$ into simpler $L$-functions. For example, a still open conjecture of Dedekind asserts that $\zeta(s)$ always divides $\zeta_K(s)$, in the sense that $\zeta_K(s)/\zeta(s)$ is entire; Dedekind's conjecture is known to hold for Galois extensions. For example, if $K$ is the cyclotomic extension generated by the $q$-th roots of unity, then
\begin{equation}
\label{1-5}
\zeta_K(s) = \prod_{\chi} L(s,\chi^*),
\end{equation}
where the product is over all Dirichlet characters $\chi$ (mod $q$) and $\chi^*$ is the primitive character inducing $\chi$. The {\it Dirichlet $L$-functions} $L(s,\chi)$ are defined by
\[
L(s,\chi) = \sum_{n=1}^\infty \frac{\chi(n)}{n^s}
\]
and form, as $q$ varies over all natural numbers, another very interesting class of arithmetic $L$-functions. In some sense, such $L$-functions represent the simplest examples after $\zeta(s)$ in the hierarchy of $L$-functions, and were introduced by Dirichlet twenty years before Riemann's $\zeta(s)$, but only as functions of the real variable $s$. They also satisfy, especially when $\chi$ is primitive, analytic properties similar to those of $\zeta(s)$, and have several applications, in particular to the study of the distribution of primes in arithmetic progressions and of the class number of quadratic fields. Factorizations of type \eqref{1-5}, apart from their intrinsic interest, have a relevant role in the study of the reciprocity laws, which deal, roughly speaking, with the decomposition of primes in number fields.

\medskip
All the arithmetic $L$-functions discussed so far have been further generalized, and new classes of such $L$-functions were introduced, notably by E.Artin and A.Weil, each with specific aims and properties. All of them share, at least conjecturally, the same basic analytic structure seen above. We refer the reader to the classical books by Lang \cite{Lan/1994}, Narkiewicz \cite{Nar/2004} and Neukirch \cite{Neu/1999} for an introduction to various $L$-functions associated with algebraic number fields.

\bigskip
{\bf 1.3. Geometric $L$-functions.} Let $E$ be an elliptic curve over $\QQ$ with conductor $N$. Then $E$ has good reduction at all primes $p\nmid N$, hence $E$ becomes an elliptic curve over the finite field $\FF_p$ for such primes. For every prime $p$, let $|E/\FF_p|$ denote the number of points of $E$ over $\FF_p$ and let
\[
a_p = p+1-|E/\FF_p|.
\]
For $p\nmid N$ the {\it local zeta function} of $E/\FF_p$ is defined as
\[
Z(X,E/\FF_p) = \frac{1 - a_pX+pX^2}{(1-X)(1-pX)},
\]
and in the 1930's Hasse proved the bound $|a_p|\leq 2p^{1/2}$. Writing $a(p)=a_pp^{-1/2}$, this implies that
\[
1 - a(p)X+X^2 = (1-\alpha_p X)(1- \overline{\alpha_p} X) \hskip1.5cm |\alpha_p| = 1.
\]
Moreover, for $p|N$ one has $|a(p)|\leq 1$. The (normalized) {\it Hasse-Weil $L$-function} of $E$ is defined as
\begin{equation}
\label{RHell3}
L(s,E) = \prod_{p|N}\big(1-\frac{a(p)}{p^s}\big)^{-1} \prod_{p\nmid N}\big(1-\frac{\alpha_p}{p^s}\big)^{-1} \big(1-\frac{\overline{\alpha_p}}{p^s}\big)^{-1},
\end{equation}
the Euler product being absolutely convergent for $\sigma > 1$. 

\medskip
The ideas and techniques used to get analytic continuation and functional equation for the arithmetic $L$-functions, based on the existence of certain theta-functions counterparts of such $L$-functions, do not apply to the Hasse-Weil $L$-functions. The analytic properties were obtained in few special cases by showing that $L(s,E)$ coincides with another $L$-function, of different nature, for which such properties were already known. In the 1950's this approach has been codified by the {\it Shimura-Taniyama conjecture}, asserting that the Hasse-Weil $L$-function of any elliptic curve coincides with a suitable modular $L$-function. We shall meet these $L$-functions in a moment, but here we anticipate that they possess the required analytic properties. In 1995 A.Wiles, partly in collaboration with R.Taylor, gave a proof of the famous Fermat Last Theorem along the lines first suggested by G.Frey, by proving the Shimura-Taniyama conjecture in a particularly important special case. The full proof of the conjecture was obtained in 2001 by C.Breuil, B.Conrad, F.Diamonf and R.Taylor, still relying on Wiles' ideas. As a consequence, it is now known that the Hasse-Weil $L$-functions $L(s,E)$ are entire and satisfy the functional equation
\[
\big(\frac{\sqrt{N}}{2\pi}\big)^s \Gamma(s+\frac12) L(s,E) = \pm \big(\frac{\sqrt{N}}{2\pi}\big)^{1-s} \Gamma(1-s+\frac12) L(1-s,E).
\]

\medskip
The proof of the Shimura-Taniyama conjecture is of great significance for another important conjecture, dealing with the arithmetic of elliptic curves. It was shown by L.J.Mordell in the 1920's that the group $E(\QQ)$ of the rational points on $E$ is finitely generated. The torsion part of $E(\QQ)$ is rather well understood at present, thanks to the work of B.Mazur in the 1970's, but the determination of the rank $r(E)$ of $E(\QQ)$ is still an open problem. The famous {\it Birch and Swinnerton-Dyer conjecture} proposes a deep link between $r(E)$ and $L(s,E)$, namely that $r(E)$ coincides with the order of vanishing of $L(s,E)$ at $s=1/2$, the central point of the critical strip (recall that $L(s,E)$ is normalized). The conjecture has been verified in certain cases, and may be viewed as an instance of the local-global principle in arithmetic geometry. Indeed, the Hasse-Weil $L$-function is defined by means of local data and carries global information on the arithmetic of the elliptic curve. We conclude by mentioning that a recent breakthrough on the rank of elliptic curves has been obtained by M.Bhargava and his collaborators, by a different approach.

\medskip
The Hasse-Weil $L$-function has been defined, in an analogous way, for any smooth projective variety. Little is known unconditionally about such functions, but there are standard conjectures describing their analytic continuation and functional equation. Moreover, we shall see in the next subsection that there is an analog of the Shimura-Taniyama conjecture for all Hasse-Weil $L$-functions. We refer the reader to the elementary introduction to elliptic functions and their $L$-functions by Lozano-Robledo \cite{Lo-Ro/2011}, to Wiles \cite{Wil/2006} for a survey of the Birch and Swinnerton-Dyer conjecture, to Bertolini-Canuto \cite{Be-Ca/1996} for a survey of the proof of the Shimura-Taniyama conjecture and to Serre \cite{Ser/1970} for the standard conjectures in the algebraic varieties case.

\bigskip
{\bf 1.4. Modular $L$-functions.} A different type of $L$-functions, associated with modular forms, was introduced by E.Hecke in the 1930's. Here we consider the simplest case of $L$-functions associated with cusp forms of integral weight $k$ for the modular group $\Gamma=$ SL$(2,\ZZ)$. A cusp form $f(z)$ of weight $k$ for $\Gamma$ is holomorphic on the upper half-plane $\H=\{z=x+iy: y>0\}$, satisfies
\[
f\big(\frac{az+b}{cz+d}\big) = (cz+d)^kf(z) \quad \text{for every} \  \left( {a\atop c} \ {b\atop d} \right) \in \Gamma
\]
and is holomorphic and vanishing at $\infty$. In particular, $f(z)$ has a Fourier expansion of type
\begin{equation}
\label{1-new2}
f(z) = \sum_{n=1}^\infty a_ne^{2\pi inz}.
\end{equation}
The cusp forms of weight $k$ for $\Gamma$ form a finite dimensional linear space over $\CC$, denoted by $S_k(\Gamma)$. Its dimension can be explicitly computed; in particular, $S_k(\Gamma)$ is trivial unless $k$ is an even integer $\geq 12$. A classical example is $\Delta(z)$, of weight 12, whose Fourier coefficients are given by the famous Ramanujan $\tau$-function. Writing $a(n) = a_nn^{-(k-1)/2}$, the (normalized) {\it Hecke $L$-function} associated with $f\in S_k(\Gamma)$ is defined by the Dirichlet series
\begin{equation}
\label{1-new3}
L(s,f) = \sum_{n=1}^\infty \frac{a(n)}{n^s},
\end{equation}
which turns out to be absolutely convergent for $\si>1$.

\medskip
The analytic properties of $L(s,f)$ were obtained by Hecke by means of a simple extension of Riemann's method: $L(s,f)$ is entire and satisfies the functional equation
\begin{equation}
\label{1-new4}
\Lambda(s,f) = i^k\Lambda(1-s,f), \quad \text{where} \ \Lambda(s,f) = \big(\frac{1}{2\pi}\big)^s\Gamma\big(s+\frac{k-1}{2}\big) L(s,f).
\end{equation}
In general, the Hecke $L$-functions $L(s,f)$ don't have an Euler product. However, Hecke gave an interesting characterization of the cusp forms whose $L$-function has this property by means of the Hecke operators $T_n$, $n=1,2,\dots$, which send $f\in S_k(\Gamma)$ to $T_nf\in S_k(\Gamma)$. Precisely, $L(s,f)$ has Euler product if and only if $f(z)$ is a normalized (i.e. $a_1=1$) simultaneous eigenfunction of all $T_n$. Moreover, $S_k(\Gamma)$ has a unique basis of normalized eigenfunctions for all $T_n$. Hecke's theory detects also the shape of the Euler factors, namely
\begin{equation}
\label{1-6}
L(s,f) = \prod_p\big(1 - \frac{a(p)}{p^s} + \frac{1}{p^{2s}}\big)^{-1}.
\end{equation}
The finer shape of \eqref{1-6} was predicted by the Ramanujan conjecture and proved by P.Deligne in 1974, as a consequence of the Riemann Hypothesis for varieties over finite fields; precisely, \eqref{1-6} can be written as
\begin{equation}
\label{1-7}
L(s,f) = \prod_p\big(1 -\frac{\alpha_p}{p^s}\big)^{-1} \big(1 -\frac{\overline{\alpha_p}}{p^s}\big)^{-1} \hskip1.5cm |\alpha_p|=1.
\end{equation}
Note that \eqref{RHell3}, which precedes Deligne's proof, is similar to \eqref{1-7} when $k=2$ and $N=1$. Moreover, the above results show that the Hecke $L$-functions, particularly those associated with eigenforms, share the same basic analytic properties of the arithmetic and geometric $L$-functions. This may appear somewhat unexpected, since the nature of these classes of $L$-functions looks quite different, but the Shimura-Taniyama conjecture shows that, in certain cases, there are deep reasons for such similarities; we shall say more on this in a moment.

\medskip
Of course, due to the above similarities, analogous problems and applications may be formulated for the Hecke $L$-functions. Moreover, the finer theory of modular forms, including the non-holomorphic case, has deep connections with several arithmetical problems. One side of this connections is the link between modular forms and certain exponential sums which play an important role in number theory. 

\medskip
Another interesting aspect of the $L$-functions associated with modular forms, perhaps less visible in the context of the $L$-functions treated in the previous subsections, is related to certain structural properties. For example, given two eigenforms $f,g\in S_k(\Gamma)$ one may form the {\it Rankin-Selberg convolution $L$-function} $L(s,f\otimes \overline{g})$. Denoting by $\alpha_i(p)$, $i=1,2$, the two $\alpha$-coefficients in \eqref{1-7}, and defining analogously $\beta_j(p)$, $j=1,2$, for $g(z)$, such an $L$-function is then defined as
\begin{equation}
\label{1-8}
L(s,f\otimes \overline{g}) = \prod_p \prod_{i,j=1}^2 \big(1 -\frac{\alpha_i(p)\overline{\beta_j(p)}}{p^s}\big)^{-1}.
\end{equation}
By a brilliant two-dimensional extension of Riemann's method, R.A.Rankin and A.Selberg, independently and roughly at the same time (around 1940), were able to obtain the meromorphic continuation and functional equation for $L(s,f\otimes \overline{g})$. Once again, the $L$-function $L(s,f\otimes \overline{g})$ has at most a simple pole at $s=1$, and the functional equation relates its values at $s$ and at $1-s$ via a suitable product of $\Gamma$-factors. It is interesting to note that the polar structure of $L(s,f\otimes \overline{g})$ reflects a kind of orthogonality among the eigenforms of $S_k(\Gamma)$, as well as among their $L$-functions: the residue at $s=1$ of $L(s,f\otimes \overline{g})$ does not vanish if and only if $f=g$. Such an orthogonality, which later turned out to be a rather general phenomenon among $L$-functions, has deep consequences in various aspects of the theory. Due to the shape of \eqref{1-7}, when $f=g$ the product over $i$ and $j$ in \eqref{1-8} contains two factors of Riemann $\zeta$-type, see \eqref{1-1}. Therefore, one may wonder if $L(s,f\otimes \overline{f})$ is divisible by $\zeta(s)$ (or even by $\zeta(s)^2$). The answer was given by G.Shimura in 1973, who showed that
\begin{equation}
\label{1-9}
L(s,f\otimes \overline{f}) = \zeta(s) L(s,\text{sym}^2f),
\end{equation}
where the {\it symmetric square $L$-function} $L(s,\text{sym}^2f)$, implicitly defined by \eqref{1-8} and \eqref{1-9}, has analytic continuation over $\CC$ and satisfies a functional equation (but is not further divisible by $\zeta(s)$). Higher order convolutions and symmetric powers can be considered, but here our knowledge is not yet complete. A further structural property, still of convolution nature, is given by the {\it twisted Hecke $L$-function}
\begin{equation}
\label{1-new5}
L(s,f\otimes\chi) = \sum_{n=1}^\infty \frac{a(n)\chi(n)}{n^s},
\end{equation}
where $\chi$ is a primitive Dirichlet character. Again, $L(s,f\otimes\chi)$ has analytic continuation and functional equation, and of course has the Euler product if $L(s,f)$ does.

\medskip
The above $f\otimes\overline{g}$, sym$^2f$ and $f\otimes\chi$ notation is reminiscent of representation theory. Actually, a vast reinterpretation and extension of the theory of modular forms has been developed starting from the 1960's, in particular by R.P.Langlands, leading the the theory of automorphic representations, automorphic forms and their associated {\it automorphic $L$-functions}. This is a difficult and highly technical area, but the output is rewarding. In particular, we have now a very general class of $L$-functions satisfying standard analytic properties, and the above convolution operations find a natural interpretation here; the Hecke $L$-functions associated with eigenforms are a very special case. However, the most interesting side of this theory is the impressive unifying role predicted by the {\it Langlands program}. We already reported that the Shimura-Taniyama conjecture, which is now a theorem, asserts that every elliptic $L$-function coincides with a suitable modular $L$-function. In much greater generality, the Langlands program predicts that every arithmetic or geometric $L$-function is in fact an automorphic $L$-function. In other words, according to the Langlands program such $L$-functions, a priori of very different nature, can all be realized as $L$-functions associated with suitable automorphic representations. The implications of such a visionary conjecture are clearly of fundamental importance. The Langlands program has been verified in few cases, the Shimura-Taniyama conjecture being at present the most spectacular.

\medskip
The interested reader is referred to the books by Hecke \cite{Hec/1983} and Iwaniec \cite{Iwa/1997} for the theory of modular forms and their $L$-functions. Also, we refer to the collection of articles Bernstein-Gelbart \cite{Be-Ge/2004} and to Bombieri \cite{Bom/2005} for an introduction to the Langlands program.

\bigskip
We conclude pointing out that $L$-functions have several amazing properties. To quote one, concerning the Riemann zeta function but shared by many other $L$-functions, we mention the {\it universality} of $\zeta(s)$, discovered by S.M.Voronin in 1975. Let $f(s)$ be holomorphic and non-vanishing on a closed disk $D$ inside the strip $1/2<\si<1$, and let $\epsilon>0$. Then there exists a real number $\tau$ (actually, there are infinitely many such $\tau$'s) with the property that 
\[
\max_{s\in D}|\zeta(s+i\tau) - f(s)| <\epsilon.
\]
In other words, one can approximate with any precision any holomorphic non-vanishing function on any closed disk as above by means of vertical shifts of $\zeta(s)$.

\bigskip
\section{Classical converse theorems}

\smallskip
We have seen that $L$-functions can be associated with several different kinds of structures; the aim of converse theorems is to detect the nature of such structures starting from the analytic properties of the $L$-functions. Clearly, there are intrinsic limitations to the effectiveness of converse theorems; for example, we already pointed out that there exist non-isomorphic number fields with the same Dedekind zeta function. Nevertheless, converse theorems often provide interesting characterizations, which turn out to be useful in several frameworks.

\bigskip
{\bf 2.1. Hamburger's theorem.} It is not surprising that the first converse theorem was obtained in connection with the Riemann zeta function. Indeed, in 1921 H.Hamburger obtained (a slightly more general version of) the following result. 

\medskip
{\bf Theorem 2.1.} (Hamburger \cite{Ham/1921}) {\sl Let $F(s)$ be a Dirichlet series, absolutely convergent for $\si>1$, such that $(s-1)^mF(s)$ is entire of finite order for some integer $m$. Moreover, suppose that $F(s)$ satisfies the functional equation
\begin{equation}
\label{2-1}
\pi^{-s/2} \Gamma\big(\frac{s}{2}\big)F(s) = \pi^{-(1-s)/2} \Gamma\big(\frac{1-s}{2}\big)F(1-s).
\end{equation}
Then $F(s) = c\zeta(s)$ for some $c\in\CC$.} 

\medskip
Thus, in view of \eqref{1-new1}, Hamburger's theorem essentially asserts that the Riemann zeta function is characterized by its functional equation inside a rather large class of Dirichlet series. There are currently several proofs of Hamburger's theorem, based on few principles. Some proofs are based on an analogous uniqueness statement for certain theta functions, other proofs exploit certain periodicity properties induced by the special form of \eqref{2-1}. We refer to Chapter II of Titchmarsh \cite{Tit/1986} and to Piatetski-Shapiro and Raghunathan \cite{PS-Ra/1995} for a discussion of Hamburger's theorem, and to Section 3 for a more recent approach in a general framework, giving Theorem 2.1 as a very special case.

\medskip
The key assumption in Hamburger's theorem is the functional equation, a distinctive property of $\zeta(s)$ and, in general, of all $L$-functions; most classical converse theorems are based on assumptions of this type. Clearly, the solutions $F(s)$, with suitable analytic requirements as in Theorem 2.1, of a general functional equation of Riemann type (see various examples in Section 1) form a linear space $W$, say; hence characterizing uniquely a given $L$-function by means of such data is simply impossible. We may therefore regard Hamburger's theorem as a perfect converse theorem, since $\zeta(s)$ is characterized with maximal precision. 

\medskip
Actually, Theorem 2.1 says that the dimension of the space $W$ is 1 in this case. This is a very rare event, in the sense that usually such spaces $W$ have larger dimension, often infinite, in the (few !) known cases; see the next subsection. Since $L$-functions come usually from finite dimensional structures, we see that converse theorems must involve non-trivial information. Moreover, the dimension of $W$ is highly sensitive to variations in the data of the functional equation; for example, slight modifications to \eqref{2-1} give already a different output. Indeed, write $\overline{F}(s) = \overline{F(\overline{s})}$ and consider the space $W$ of solutions of the functional equation
\begin{equation}
\label{2-2}
\big(\frac{q}{\pi}\big)^{s/2} \Gamma\big(\frac{s+a}{2}\big)F(s) = \omega \big(\frac{q}{\pi}\big)^{(1-s)/2} \Gamma\big(\frac{1-s+a}{2}\big)\overline{F}(1-s).
\end{equation}
For $q\in\NN$, $a=0$ or 1 and suitable choices of $\omega$ of modulus 1, \eqref{2-2} is satisfied by the Dirichlet $L$-functions $L(s,\chi)$ with primitive $\chi$ (mod $q$), and $W$ has in these cases dimension which grows linearly with $q$. Again we refer to \cite{PS-Ra/1995} for a discussion of this issue. 

\medskip
We conclude observing that if one adds in Theorem 2.1 the further requirement that $F(s)$ is an Euler product, the other distinctive property of $L$-functions, then it follows immediately that $F(s)=\zeta(s)$. The situation is different for \eqref{2-2}, where functions $L(s,\chi)$ with distinct $\chi$'s may satisfy the same functional equation; Kaczorowski-Molteni-Perelli \cite{K-M-P/2010} characterized the values of $q$ for which this happens, and Molteni \cite{Mol/2010},\cite{Mol/2010a} gave quantitative multiplicity results.

\bigskip
{\bf 2.2. Hecke's converse theorem.} Hamburger's theorem remained a somewhat isolated result until Hecke reconsidered it from the point of view of his theory of modular forms and associated $L$-functions. In Section 1 we briefly recalled Hecke's theory for the modular group $\Gamma$, in the case of cusp forms; we still continue with this framework for ease of exposition, since the slightly more general case of modular forms is dealt with similarly. In particular, we recall that given $f\in S_k(\Gamma)$ with Fourier series as in \eqref{1-new2}, the (normalized) Hecke $L$-function \eqref{1-new3} is entire of finite order and satisfies the functional equation \eqref{1-new4}. Hecke proved that the opposite implication holds true as well, thus showing a perfect correspondence between cusp forms for $\Gamma$ and Dirichlet series $F(s)$ with analytic continuation satisfying \eqref{1-new4}. Actually, Hecke proved a more general correspondence which led him, for example, to an alternative approach to Hamburger's theorem. Given a positive even integer $k$ and a sequence of complex numbers $a(n)$, we write $a_n=a(n)n^{(k-1)/2}$ and form the two series
\begin{equation}
\label{2-new1}
F(s) = \sum_{n=1}^\infty \frac{a(n)}{n^s} \qquad \text{and} \qquad f(z) = \sum_{n=1}^\infty a_n e^{2\pi inz},
\end{equation}
and consider
\[
\Lambda(s,F) = \big(\frac{1}{2\pi}\big)^s\Gamma\big(s+\frac{k-1}{2}\big) F(s).
\]
With this notation we have

\medskip
{\bf Theorem 2.2.} (Hecke's converse theorem) {\sl Suppose that $F(s)$ is absolutely convergent for $\si>1$ and has continuation over $\CC$ as an entire function of finite order. Moreover, suppose that $F(s)$ satisfies the functional equation
\begin{equation}
\label{2-new2}
\Lambda(s,F) = i^k\Lambda(1-s,F).
\end{equation}
Then $f\in S_k(\Gamma)$.}

\medskip
Theorem 2.2 can be proved by a standard technique, extending one of Riemann's proofs of the functional equation of $\zeta(s)$, since $F(s)$ and $f(z)$ are linked by the Mellin transform; nevertheless it is an interesting result. First note that it is less precise than Hamburger's theorem: from the above analytic properties we can only say, in general, that $F(s)$ is the $L$-function of {\it some} $f\in S_k(\Gamma)$. But we already noticed that this is an intrinsic limitation since the solutions of \eqref{2-new2} form a linear space $W$, and the limitation is stronger when $W$ has higher dimension. The interesting information given by Theorem 2.2 is that, essentially, the modularity of $f(z)$ is {\it equivalent} to the functional equation of $L(s,f)$. Moreover, such an equivalence allows to compute the dimension of $W$, by looking at the dimension of its counterpart $S_k(\Gamma)$. This possibility is very interesting since, as far as we know, at present there is no direct approach to the dimension of $W$ from the Dirichlet series side, without involving in some way the modularity properties of the functions $f(z)$.

\medskip
The reason for such an equivalence is quite simple. The modular group $\Gamma$ has two generators, which in terms of their action on the upper half-plane $\H$ are described as the shift $z\mapsto z+1$ and the inversion $z\mapsto -1/z$; to check the modularity of $f(z)$ for $\Gamma$, it is clearly enough to check its modularity on the generators. Now, roughly speaking, the shift-modularity is equivalent to the Fourier expansion in \eqref{1-new2}, while the inversion-modularity is equivalent to functional equation \eqref{1-new4}. Actually, this is a simple but important special case of a duality principle, first studied by S.Bochner in the 1950's, between Dirichlet series with a general Riemann type functional equation and their Mellin transforms; see, e.g., Kaczorowski-Perelli \cite{Ka-Pe/2011b}. Another instance of such a duality is represented by Hecke's correspondence, already mentioned above, between modular forms for the triangle groups $G(\lambda)$, $\lambda\in\RR$, generated by $z\mapsto z+\lambda$ and $z\mapsto -1/z$, and Dirichlet series with a functional equation of type similar to \eqref{1-new4}; the case $\lambda=1$ corresponds to the modular forms for $\Gamma$. The more interesting part of such a correspondence is that Hecke managed to compute exactly the dimension of the resulting spaces $W$, again arguing on the modular forms side. It turns out that for $0<\lambda \leq 2$ the dimension of $W$ is finite (usually 0), while for $\lambda>2$ all $W$'s are infinite-dimensional; moreover, there is a natural reason for this phenomenon, depending on the geometry of the fundamental domain of $G(\lambda)$. Hecke's theory is presented in detail in Hecke \cite{Hec/1983} and Ogg \cite{Ogg/1969}.

\bigskip
{\bf 2.3. Weil's converse theorem.} All aspects of Hecke's theory become more complicated passing from the modular group $\Gamma$ to the Hecke congruence subgroups of level $N$, denoted by $\Gamma_0(N)$ and consisting of the matrices $(\begin{smallmatrix}a&b\\c&d\end{smallmatrix})\in\Gamma$ with $c$ divisible by $N$; note that $\Gamma_0(1)=\Gamma$. Such subgroups are of arithmetical interest, and a suitable extension to $\Gamma_0(N)$ of Hecke's theory for $\Gamma$ exists. Very briefly, we recall that for integer weight $k$ one is led to consider spaces of modular or cusp forms with multiplier $\psi$, the latter being denoted by $S_k(\Gamma_0(N),\psi)$, where $\psi$ is a Dirichlet character (mod $N$). Moreover, each cusp form has a Fourier expansion as in \eqref{1-new2} and an associated (normalized) Hecke $L$-function as in \eqref{1-new3}, with analytic continuation and functional equation. As in the case of $\Gamma$, the spaces $S_k(\Gamma_0(N),\psi)$ are always finite dimensional. However, there are some important differences; for simplicity we often consider the case where $\psi$ is trivial, and denote the corresponding space by $S_k(\Gamma_0(N))$.

\medskip
When $f\in S_k(\Gamma_0(N))$, the associated $L(s,f)$ is entire of finite order and the shape of the functional equation is
\begin{equation}
\label{2-3}
\Lambda(s,f) = i^k\Lambda(1-s,f_{|\omega_N}), \quad \text{where} \ \  \Lambda(s,f) = \big(\frac{\sqrt{N}}{2\pi}\big)^s\Gamma\big(s+\frac{k-1}{2}\big) L(s,f)
\end{equation}
and $f_{|\omega_N}(z)$ is another cusp form, also belonging to $S_k(\Gamma_0(N))$. More precisely, $\omega_N$ is the Fricke involution and 
\[
f_{|\omega_N}(z) = (\sqrt{N}z)^{-k} f(-1/Nz).
\]
We note in passing that when $f(z)$ is a newform of $S_k(\Gamma_0(N),\psi)$ we have
\[
f_{|\omega_N}(z) = \epsilon(f) \overline{f}(z),
\]
where $|\epsilon(f)|=1$ and $\overline{f}(z)$ denotes the series in \eqref{1-new2} with conjugate coefficients; newforms are the analog for $\Gamma_0(N)$ of the eigenforms for $\Gamma$, and $\overline{f}\in S_k(\Gamma_0(N),\overline{\psi})$. In this case $L(s,f)$ has an Euler product, and functional equation \eqref{2-3} takes the nicer form
\[
\Lambda(s,f) = \omega \Lambda(1-s,\overline{f})
\]
where $|\omega|=1$. This observation will be relevant in Section 3. 

\medskip
For levels  $N=1, 2, 3$ or 4 we have the exact analog of Theorem 2.2, namely the functional equation of $L(s,f)$ is enough to characterize the cusp forms in $S_k(\Gamma_0(N))$. This is due to the fact that in these cases $\Gamma_0(N)$ has two generators, which are very close to those of $\Gamma$; but for $N\geq 5$ the situation changes dramatically. Indeed, on the one hand the groups $\Gamma_0(N)$ have in general many generators, and one functional equation is not anymore enough to ensure modularity for all generators. On the other hand, the spaces $W$ associated with \eqref{2-3} are infinite-dimensional, since they correspond to the case $\lambda=\sqrt{N}>2$ in Hecke's theory for the groups $G(\lambda)$. As a consequence, roughly speaking, there are more Dirichlet series with functional equation than cusp forms. Hence something more is needed for their characterization.

\medskip
Another interesting part of Hecke's theory concerns the {\it twist} of modular forms and their $L$-functions. Given $f\in S_k(\Gamma_0(N),\psi)$ and a primitive character $\chi$ (mod $m$) with $(m,N)=1$, Hecke considered the twisted cusp form
\[
f^\chi(z) = \sum_{n=1}^\infty a_n\chi(n) e^{2\pi inz}
\]
and the associated twisted $L$-function $L(s,f\otimes\chi)$ defined by \eqref{1-new5}. Then, by an argument based on Gauss sums, he showed that $f^\chi\in S_k(\Gamma_0(Nm^2),\psi\chi^2)$. Therefore, thanks to his general theory, the twisted $L$-function is entire of finite order and satisfies the functional equation
\begin{equation}
\label{2-4}
\Lambda(s,f\otimes\chi) = \omega(\chi) \Lambda(1-s,f_{|\omega_N}\otimes\overline{\chi}),
\end{equation}
where
\[
\Lambda(s,f\otimes\chi) = \big(\frac{\sqrt{N}m}{2\pi}\big)^s\Gamma\big(s+\frac{k-1}{2}\big) L(s,f\otimes\chi)
\]
and $\omega(\chi)$ is a certain complex number of modulus 1 depending on $\psi$, $\chi$ and $k$.

\medskip
In 1967, A.Weil had the idea of using twists in order to characterize higher level modular forms, proving his famous converse theorem; again, for simplicity we state it in a slightly specialized form, which however presents all the main features of the general case. Let $N$ and $k$ be positive integers, $\psi$ be a character (mod $N$), $a(n)$ and $b(n)$ be sequences of complex numbers and write $a_n=a(n)n^{(k-1)/2}$, $b_n=b(n)n^{(k-1)/2}$. Moreover, consider the series
\[
F(s) = \sum_{n=1}^\infty \frac{a(n)}{n^s}, \ G(s) = \sum_{n=1}^\infty \frac{b(n)}{n^s}, \
 f(z) = \sum_{n=1}^\infty a_n e^{2\pi inz} \ \text{and} \ g(z) = \sum_{n=1}^\infty b_n e^{2\pi inz},
\]
and the function
\[
\Lambda(s,F) = \big(\frac{\sqrt{N}}{2\pi}\big)^s\Gamma\big(s+\frac{k-1}{2}\big) F(s),
\]
and analogously for $G(s)$. Consider also the twisted functions
\[
F(s,\chi) = \sum_{n=1}^\infty \frac{a(n)\chi(n)}{n^s} \quad \text{and} \quad \Lambda(s,F\otimes\chi) = \big(\frac{\sqrt{N}m}{2\pi}\big)^s\Gamma\big(s+\frac{k-1}{2}\big) F(s,\chi),
\]
and analogously for $G(s)$. Further, let $\omega(\chi)$ be as in \eqref{2-4} and let $\M$ be a set of primes, all coprime with $N$, which meets every reduced residue class, i.e. for every $n\geq 1$ and $(a,n)=1$ there exists a prime $m\in\M$ with $m\equiv a$ (mod $n$). With this notation we have

\medskip
{\bf Theorem 2.3.} (Weil's converse theorem) {\sl Suppose that $F(s)$ and $G(s)$ are absolutely convergent for $\si>1$, have continuation as entire functions of finite order and satisfy
\[
\Lambda(s,F) = i^k \Lambda(1-s,G).
\]
Moreover, suppose that for every primitive character $\chi$} (mod $m$), $m\in\M$, {\sl the functions $F(s,\chi)$ and $G(s,\chi)$ have continuation as entire functions of finite order and satisfy
\[
\Lambda(s,F\otimes\chi) = \omega(\chi) \Lambda(1-s,G\otimes\overline{\chi}).
\]
Then $f\in S_k(\Gamma_0(N),\psi)$ and $g(z) = f_{|\omega_N}(z)$.}

\medskip
Thus, roughly speaking, the functional equations of $L(s,f)$ and of many of its twists characterize the cusp form $f(z)$. Weil's basic idea, whose implementation involves some delicate points, is to use the twisted functional equations to check modularity for all generators of $\Gamma_0(N)$. It is interesting that the good behavior under twists, which is one of the structural properties of modular forms, is actually sufficient for their characterization.

\medskip
Finally, note that Weil's converse theorem needs infinitely many functional equations to check modularity for finitely many generators; it is natural to ask if finitely many such functional equations are enough. The answer is affirmative, and results of this sort were later proved by M.J.Razar and W.C.W.Li. We refer to Ogg \cite{Ogg/1969} and Iwaniec \cite{Iwa/1997} for a presentation of Weil's theorem.

\bigskip
{\bf 2.4. Other converse theorems.} In Section 1 we mentioned that the theory of modular forms and of their $L$-functions is a special case of the general theory of automorphic forms and $L$-functions. The question of how to recognize automorphic forms from an analytic point of view may be asked in such a general framework as well, and indeed general converse theorems exist. We already mentioned that twisting is essentially a convolution type operation; it turns out that convolutions of Rankin-Selberg type are the right replacement of twists in the general case. We do not enter detailed statements here, as we did not even give the definition of automorphic $L$-functions in Section 1. However, roughly speaking we may say that in order to characterize an automorphic form of a given degree, one must show good properties of the convolutions of the associated $L$-function with the automorphic forms of lower degree. We refer to Cogdell \cite{Cog/2007} for a survey of converse theorems for automorphic $L$-functions.

\medskip
Up to now we discussed converse theorems where one tries to characterize modular forms by means of functional equations. We already remarked that this is due to the fact that the functional equation is a distinctive property of $L$-functions; but it is not the only one. Actually, arithmetic is brought into $L$-functions mainly by the Euler product, hence one may ask if converse theorems exist based on Euler products instead of functional equations. For example, a natural question to start with is as follows. Suppose that 
\[
F(s) = \prod_{p}\big(1-\frac{\alpha_p}{p^s}\big)^{-1},
\]
with $|\alpha_p|=1$ for all primes $p$, is such that  $(s-1)^mF(s)$ is entire of finite order for some $m\in\NN$. What can we say about $F(s)$ ? Certainly $\zeta(s)$ is a solution of this problem; moreover, since $\alpha_p=e^{i\theta_p}$ with $0\leq \theta_p<2\pi$, no modification of a finite number of factors of its Euler product gives another solution. As far as we know this problem has not been investigated in the literature, but the situation changes if one is ready to consider the joint effect of Euler product and functional equation.

\medskip
In 1995, Conrey-Farmer \cite{Co-Fa/1995} proposed the study of an alternative version of Weil's converse theorem, where the properties of the twisted $L$-functions are replaced by the assumption of a suitable Euler product. In our normalized setting, given positive integers $N$ and $k$ one says that the Dirichlet series $F(s)$ in \eqref{2-new1} satisfies a functional equation of {\it degree $2$, level} $N$ and {\it weight} $k$ if
\begin{equation}
\label{2-5}
\Lambda(s) = \pm(-1)^{k/2}\Lambda(1-s), \quad \text{where} \ \Lambda(s) = \big(\frac{\sqrt{N}}{2\pi}\big)^s \Gamma\big(s+\frac{k-1}{2}\big) F(s),
\end{equation}
while it has an Euler product of {\it degree $2$, level} $N$ and {\it weight} $k$ if
\[
F(s) = \prod_p F_p(s)
\]
where
\begin{equation}
\label{2-6}
\begin{split}
F_p(s) &= (1-a(p)p^{-s} +p^{-2s})^{-1} \quad  \text{if} \ p\nmid N, \\
F_p(s)&= (1- q^{-s})^{-1}  \quad  \text{if} \ p||N, \\
F_p(s)&= 1  \quad  \text{if} \ p^2|N.
\end{split}
\end{equation}
Moreover, let $f(z)$ be defined as in \eqref{2-new1}. With this notation, and once again in a slightly weaker form, the following partial result was proved in \cite{Co-Fa/1995}.

\medskip
{\bf Theorem 2.4.} (Conrey-Farmer's converse theorem) {\sl Let $5\leq N\leq 12$ or $14\leq N\leq 17$ or $N=23$. Suppose that $F(s)$ is absolutely convergent for $\si>1$ and has continuation as an entire function of finite order. Suppose that $F(s)$ has both a functional equation and an Euler product of degree $2$, level $N$ and weight $k$. Then $f\in S_k(\Gamma_0(N))$.}

\medskip
Hence the twists in Weil's theorem are replaced, for the above special values of $N$, by an Euler product of the expected shape. A further case, namely $N=13$, has been treated by Conrey-Farmer-Odgers-Snaith \cite{C-F-O-S/2007}. Although the technique is quite different, the very basic idea in Theorem 2.4 is similar to Weil's theorem. Indeed, one tries to check modularity for all generators of $\Gamma_0(N)$, this time exploiting the relations coming from the Euler product instead of those coming from the twisted functional equations. The restriction to the above values of $N$ is due to the lack of a general strategy to deduce modularity relations from the Euler product information; the generators are indeed of special form in the cases under consideration. Extending Theorem 2.4 to all $N$'s is an interesting open problem. 

\medskip
Finally, Diaconu-Perelli-Zaharescu \cite{D-P-Z/2002} obtained a hybrid version of the converse theorems of Weil and Conrey-Farmer. Indeed, they  showed a general converse theorem for $S_k(\Gamma_0(N))$ assuming that $F(s)$ has both a functional equation and an Euler product of degree $2$, level $N$ and weight $k$ as in Conrey-Farmer's theorem, and moreover that $F(s,\chi)$ satisfies a functional equation for all primitive $\chi$'s modulo a {\it single} suitable prime $m$.

\medskip
We conclude the section recalling that many other converse theorems exist in the literature, many of which are concerned with the theory of Maass waves, the non-holomorphic counterpart of classical cusp forms. Some of these results are subsumed into the general converse theorems for automorphic $L$-functions, but other deal with variants of the classical cases in different frameworks.

\bigskip
\section{Converse theorems for the Selberg class}

\smallskip
In the previous sections we have seen several different kinds of $L$-functions; since the Langlands program asserts that, essentially, all classical $L$-functions belong to the class of automorphic $L$-function, one may regard the latter as the class of all $L$-functions. On the other hand, we have also seen that the classical $L$-functions share, at least conjecturally, few basic analytic properties; hence one may consider a class of $L$-functions defined axiomatically by such properties, and regard this as the class of all $L$-functions. Both points of view provide suitable environments and a wealth of interesting problems. In this section we describe some recent (and less recent) work stemming from the second approach.

\bigskip
{\bf 3.1. The Selberg class of $L$-functions.} Selberg \cite{Sel/1989} defined an axiomatic class of $L$-functions, which is now called the {\it Selberg class} and denoted by $\S$, by the following axioms. Recalling our notation $\overline{F}(s)=\overline{F(\overline{s})}$, a function $F(s)$ belongs to $\S$ if

\smallskip
{\sl (i)} ({\it Dirichlet series}) $F(s)$ is represented by a Dirichlet series
\[ 
F(s) = \sum_{n=1}^\infty \frac{a(n)}{n^s}
\]
absolutely convergent for $\sigma>1$;

{\sl (ii)} ({\it analytic continuation}) there exists an integer $m\geq 0$ such that $(s-1)^{m}F(s)$ is an entire function of finite order;

{\sl (iii)} ({\it functional equation}) $F(s)$ satisfies a functional equation of type 
\begin{equation}
\label{2-6bis}
Q^s \lprod_{j=1}^r \Gamma (\lambda_j s+ \mu_j)F(s) = \omega  Q^{1-s} \lprod_{j=1}^r \Gamma (\lambda_j (1-s)+ \overline{\mu_j})\overline{F}(1-s)
\end{equation}
where $|\omega|=1$, $r\geq 0$, $Q>0$, $\lambda_j >0$ and $\Re  \mu_j \geq 0$;

{\sl (iv)} ({\it Ramanujan conjecture}) the coefficients satisfy $a(n)\ll n^{\varepsilon}$ for every $\varepsilon > 0$;

{\sl (v)} ({\it Euler product}) $\log F(s)$ is represented by a Dirichlet series 
\[
\log F(s) = \lsum_{n=1}^{\infty}\frac{b(n)}{n^s}
\]
with $b(n)=0$ unless $n=p^m$, $m\geq 1$, and $b(n)\ll n^{\vartheta}$ for some $\vartheta<1/2$.

\smallskip
\noindent
All the examples of classical $L$-functions with Euler product encountered in Section 1, and many more, belong to $\S$, at least conjecturally. Indeed, the Ramanujan conjecture is not yet known in general for automorphic $L$-functions, while certain arithmetic $L$-functions, such as the Artin $L$-functions, are only expected to have analytic continuation of the prescribed type. Of course there are interesting cases of $L$-functions not belonging to $\S$, for example the Hecke $L$-functions of cusp forms which are not newforms, but we may safely say that the Selberg class is a suitable {\it analytic model} for the $L$-functions of arithmetic interest.  We refer to the original paper by Selberg \cite{Sel/1989}, to the early papers on the subject by Conrey-Ghosh \cite{Co-Gh/1993} and Ram Murty \cite{Mur/1994} and to the surveys of Kaczorowski-Perelli \cite{Ka-Pe/1999b}, Kaczorowski \cite{Kac/2006} and Perelli \cite{Per/2005},\cite{Per/2004},\cite{Per/2010} for the basic theory of the Selberg class. The reader may find there all the results reported below, if no reference is given.

\medskip
Some remarks on the above axioms are in order. First, the Selberg class considers $L$-functions in {\it normalized} form, where absolute convergence holds for $\si>1$, poles are allowed only at $s=1$, the functional equation reflects values according to $s\mapsto 1-s$ and the Ramanujan conjecture is expressed by $a(n)\ll n^\epsilon$. When $r=0$ in axiom {\it (iii)}, there are no $\Gamma$-factors in the functional equation; it is known that the only function in $\S$ with such a property is the constant $1$. Moreover, axiom {\it (v)} implies that every $F\in\S$ has for $\si>1$ an Euler product of type
\[
F(s) = \prod_p F_p(s), \qquad \text{where} \ \ F_p(s) =\sum_{m=0}^\infty \frac{a(p^m)}{p^{ms}}
\]
is the $p$-th Euler factor of $F(s)$. The choice of the axioms defining $\S$ is sharp; for example, the seemingly harmless condition $\vartheta<1/2$ in {\it (v)} rules out from $\S$ all non-trivial Dirichlet polynomials, which may violate the Riemann Hypothesis. Actually, Selberg conjectured that the Riemann Hypothesis holds for all $F\in\S$, and examples violating it can be produced if any of the five axioms is withdrawn. It is not difficult to derive many other basic properties of $L$-functions from the above axioms. For example, the notions of critical strip and critical line, of trivial and non-trivial zeros, the asymptotics for the number of zeros, the convexity bounds for the order of growth on vertical strips and so on, can easily be obtained. Hence the Selberg class may serve as an excellent environment if one wants to prove general theorems on $L$-functions; however, a more interesting aspect of the Selberg class is that it offers a strategic standpoint for the study of the {\it structural properties} of $L$-functions. Here are two conjectural examples, dealing essentially with the {\it independence} of $L$-functions.

\medskip
The Selberg class is a multiplicative semigroup, and a function $F\in\S$ is called {\it primitive} if it is not a product of non-trivial functions in $\S$; for example $\zeta(s)$ is known to be primitive, as well as the Dirichlet $L$-functions $L(s,\chi)$ with primitive $\chi$ and the Hecke $L$-functions of newforms. It is also known that every $F\in\S$ has a factorization into primitive functions, and it has been conjectured that {\it factorization is unique}. The unique factorization is interesting in this context since it means, roughly, that one cannot shuffle the Euler factors of distinct primitive functions to rebuild different functions still having good analytic properties. In other words, unique factorization implies a kind of rigidity of the Euler products, and hence may be regarded as a form of independence of primitive functions. In a similar direction, the {\it multiplicity one} property, known to hold in $\S$, asserts that two distinct functions cannot differ by only finitely many Euler factors.

\medskip
Another instance is provided by the {\it Selberg orthonormality conjecture}, claiming a strong form of independence among primitive functions. Precisely, given two primitive functions $F,G\in\S$ with coefficients $a(n)$ and $b(n)$, respectively, and writing $\delta_{F,G}=1$ if $F(s)=G(s)$ and $\delta_{F,G}=0$ otherwise, the orthonormality conjecture asserts that
\[
\sum_{p\leq x} \frac{a(p)\overline{b(p)}}{p} = \delta_{F,G}\log\log x + o(\log\log x).
\]
In other words, one does not expect correlation between distinct primitive functions. Clearly, Selberg's conjecture is inspired by the properties of the Rankin-Selberg convolution, which imply the orthonormality conjecture for many classical $L$-functions. Selberg's orthonormality conjecture has several interesting consequences. Some are rather direct implications, like the unique factorization conjecture or the fact that $\zeta(s)$ is the only primitive function with a pole at $s=1$. Other consequences require a deeper analysis, like Artin's conjecture, the statistical distribution of the values and the joint universality of primitive $L$-functions, and certain independence properties of their zeros. In all cases, these results show that orthonormality is a very basic and powerful structural property of $L$-functions.

\medskip
It is not difficult to see that there are choices of the {\it data} $(Q; \lambda_1,\dots,\lambda_r;\mu_1,\dots,\mu_r;\omega)$ such that functional equation \eqref{2-6bis} has no solutions in $\S$; see the last subsection. Hence another very interesting (and difficult) structural problem is to find the conditions under which \eqref{2-6bis} has solutions in $\S$, and to describe such solutions. In other words, if one is ready to accept that $\S$ is the class of all $L$-functions, this amounts to the following natural question: are all $L$-functions already known, or there exist new $L$-functions not belonging to the known families? This is nothing but a very general form of converse problem, and for its study we must first make some order inside $\S$. To this end we introduce several {\it invariants}, since the data in {\it (iii)} are not uniquely defined for a given $F\in\S$, due to the identities satisfied by the $\Gamma$ function.

\medskip
It turns out that the main invariants in $\S$ are the {\it degree} $d_F$ and the {\it conductor} $q_F$, defined by
\[
d_F=2\sum_{j=1}^r \lambda_j \hskip2cm q_F = (2\pi)^{d_F} Q^2 \prod_{j=1}^r \lambda_j^{2\lambda_j}.
\]
These definitions may look a bit cryptic, especially the conductor, but if one computes their values for the classical $L$-functions, the meaning becomes clear. For example, $\zeta(s)$ has degree and conductor 1, the Dirichlet $L(s,\chi)$ with primitive $\chi$ (mod $q$) have degree 1 and conductor $q$, the normalized Hecke $L$-functions $L(s,f)$ have degree 2 and conductor equal to the level of $f(z)$, the Dedekind zeta function of $K$ has degree equal to the degree $[K:\QQ]$ and conductor equal to the absolute value of the discriminant of $K$, and so on. So, we may regard $d_F$ and $q_F$ as a measure of the {\it analytic} and {\it arithmetic complexity} of $F(s)$, respectively; it turns out that these invariants are of great importance for the study of $\S$. However, proving converse theorems for $\S$ is not at all an easy task, since the axioms are very general. For example, a priori we only know that $d_F$ is a non-negative real number and $q_F$ is a positive real number; hence degree and conductor partition $\S$ into a continuum of still potentially large subclasses. Therefore, several other invariants have been introduced to describe the finer structure of $\S$. Here we only mention that a complete and useful set of invariants is provided by conductor, {\it root number} $\omega_F$ and the $H$-{\it invariants} $H_F(n)$, $n=0,1,\dots$; the degree appears in this set as $H_F(0)$. Moreover, there is a rather complete theory for the invariants of $\S$.

\medskip
Having introduced the invariants, we are now in position to start the description of $\S$. Recall that all known examples of $L$-functions have integer degree, and all $L$-functions of arithmetic or geometric nature are expected to be automorphic. This suggested the two main conjectures on the structure of the Selberg class. The first, called the {\it degree conjecture}, is widely accepted and easy to state; it asserts that $d_F$ is a non-negative integer for every $F\in\S$. We shall not give a precise statement for the second, less widely accepted, conjecture, as we did not define the automorphic $L$-functions. It may be called the {\it general converse conjecture} and, broadly speaking, asserts that all $F\in\S$ with integer degree are automorphic $L$-functions. One can state a more precise version of the general converse conjecture with the help of the invariants of $\S$; moreover, it is clear that the degree conjecture is also a kind of converse theorem. According to these conjectures we therefore have the following general expectation:

\smallskip
\centerline{{\it $\S$ coincides with the class of automorphic $L$-functions.}}

\smallskip
\noindent
In the remaining part of the paper we give a glimpse on the state of the art of converse theorems for the Selberg class.

\bigskip
{\bf 3.2. Nonlinear twists.} The main tool in the proof of the known cases of the degree conjecture and of the general converse conjecture are the {\it nonlinear twists} of functions in $\S$, although in few cases the original proofs were different. The theory of nonlinear twists, as well as most results presented in the next subsection, has been developed by J.Kaczorowski and the author in a series of papers \cite{Ka-Pe/1999a},\cite{Ka-Pe/2002},\cite{Ka-Pe/2005},\cite{Ka-Pe/2011a},\cite{Ka-Pe/2015},\cite{Ka-Pe/resoI},\cite{Ka-Pe/resoII},\cite{Ka-Pe/weak}. The nonlinear twists are somewhat different from the more arithmetic twists by Dirichlet characters seen before, for example in Weil's converse theorem, and indeed the ideas and techniques leading to the converse theorems in $\S$ are quite different from those used in the classical converse theorems.

\medskip
Let $F\in\S$ and write $e(x)=e^{2\pi ix}$. Given a {\it twist function} $f(n,\bfalpha)$ of type
\begin{equation}
\label{3-1}
f(n,\bfalpha) = \sum_{j=0}^N \alpha_j n^{\kappa_j} \hskip1.5cm N\geq 0, \ 0\leq\kappa_N<\dots<\kappa_0, \ \alpha_0\cdots\alpha_N\neq0,
\end{equation}
the associated nonlinear twist is defined for $\si>1$ as
\[
F(s;f) = \sum_{n=1}^\infty \frac{a(n)}{n^s} e(-f(n,\bfalpha)).
\]
The main problem about nonlinear twists is to detect their analytic properties; it turns out that our current knowledge on $F(s;f)$ depends critically on the value of the {\it leading exponent} $\kappa_0$. Very briefly, at present we understand the basic features of $F(s;f)$ whenever $\kappa_0\leq 1/d_F$, while for $\kappa_0>1/d_F$ we know certain general properties, in particular an interesting transformation formula, but only for special classes of twist functions there is a satisfactory control on the analytic properties of $F(s;f)$. The reason for such a difference depends ultimately on the different behavior of certain hypergeometric functions which arise naturally in connection with the functional equations of the $F\in\S$. Nevertheless, even such an incomplete picture is often helpful to understand the structure of $\S$.

\medskip
Two special cases of nonlinear twists turn out to be particularly interesting and useful, namely the {\it linear twist}
\[
F(s,\alpha) = \sum_{n=1}^\infty \frac{a(n)}{n^s} e(-\alpha n)
\]
and the {\it standard twist} 
\[
F^{1/d_F}(s,\alpha) = \sum_{n=1}^\infty \frac{a(n)}{n^s} e(-\alpha n^{1/d_F}).
\]
Clearly, when $F(s)$ has degree 1 we have $F^1(s,\alpha)=F(s,\alpha)$; {\it when there is no risk of confusion, we denote the standard twist simply by} $F(s,\alpha)$. Here we report only the main properties of nonlinear twists, and refer the interested reader to the above quoted papers for a more complete account. We consider only the case $\alpha_0>0$, otherwise we may simply take conjugates since $\overline{F}\in\S$ if $F\in\S$. Moreover, we always assume that $d_F>0$ since the only function in $\S$ with $d_F=0$ is the constant 1 (see Theorem 3.4). Further, for simplicity we also assume that $F(s)$ is normalized by the condition $\theta_F=0$, where the invariant $\theta_F\in\RR$ (the {\it internal shift}) controls, roughly speaking, the vertical translation of $F(s)$; the classical $L$-functions have $\theta_F=0$.

\medskip
The basic properties of the standard twist $F(s,\alpha)$ are quite well understood, and are expressed by means of the notion of {\it spectrum} of $F\in\S$, defined as
\[
\text{Spec$(F) = \{\alpha>0: a(n_\alpha)\neq 0\}$, where $n_\alpha = q_Fd_F^{-d_F}\alpha^{d_F}$ and $a(n_\alpha)=0$ if $n_\alpha\not\in\NN$.}
\]

\medskip
{\bf Theorem 3.1.} (Kaczorowski-Perelli \cite{Ka-Pe/2005},\cite{Ka-Pe/resoI}) {\sl Let $F\in\S$ with $d_F>0$ and $\theta_F=0$, and let $\alpha>0$. Then the standard twist $F(s,\alpha)$ is meromorphic over $\CC$, and is entire if $\alpha\not\in$ {\rm Spec}$(F)$. Moreover, if $\alpha\in$ {\rm Spec}$(F)$ then $F(s,\alpha)$ has at most simple poles at the points
\[
s_k = \frac{d_F+1}{2d_F} -\frac{k}{d_F} \hskip2cm k=0,1,\dots,
\]
with}
\[
\res_{s=s_0} F(s,\alpha) = c_0(F) \frac{\overline{a(n_\alpha)}}{n_\alpha^{1-s_0}} \hskip2cm c_0(F)\neq0.
\]

\medskip
It is also known that $F(s,\alpha)$ has polynomial growth on vertical lines, but here we are mainly interested in its polar structure. In particular, the pole at $s=s_0$ and the value of the residue play a relevant role in converse theorems. We remark here that, in general, {\it if $\theta_F\neq0$ then the poles of nonlinear twists, if any, are shifted by a quantity proportional to $i\theta_F$}. The properties of $F(s;f)$ in the remaining cases with leading coefficient $\kappa_0\leq 1/d_F$ are described by

\medskip
{\bf Theorem 3.2.} (Kaczorowski-Perelli \cite{Ka-Pe/resoI}) {\sl Let $F\in\S$ with $d_F>0$ and $\theta_F=0$, and let $f(n,\bfalpha)$ be as in \eqref{3-1} with $\alpha_0>0$. If $\kappa_0<1/d_F$, or $\kappa_0=1/d_F$ and $N\geq 1$, then $F(s;f)$ is entire.}

\medskip
Thus, when $\kappa_0\leq 1/d_F$, the presence of at least one term in $f(n,\bfalpha)$ with exponent $<1/d_F$ implies that the corresponding nonlinear twist has continuation to $\CC$ as an entire function. As in the case of the standard twist, we also have bounds for the order of growth on vertical lines.

\medskip
The situation is definitely more complicated when the leading exponent $\kappa_0$ is $>1/d_F$; the techniques are more difficult and the results are far less complete. The main known property in this case is a {\it transformation formula} relating a nonlinear twist $F(s;f)$ to its {\it dual} twist $\overline{F}(s;f^*)$, which we now define. Write $\omega_j=\kappa_0-\kappa_j$, $j=1,\dots,N$, and consider the semigroup
\[
\D_f = \{\omega=\sum_{j=1}^N m_j\omega_j: m_j\in\ZZ, \ m_j\geq 0\}.
\]
Moreover, let
\begin{equation}
\label{3-1bis}
\kappa_0^* = \frac{\kappa_0}{d_F\kappa_0-1}, \qquad \omega^* = \frac{\omega}{d_F\kappa_0-1}, \qquad s^* = \frac{s+\frac{d_F\kappa_0}{2}-1}{d_F\kappa_0-1}
\end{equation}
and
\[
\overline{F}(s;f^*) = \sum_{n=1}^\infty \frac{\overline{a(n)}}{n^s}e(-f^*(n,\bfalpha)), \quad \text{where} \ \
f^*(n,\bfalpha) = n^{\kappa_0^*}  \sum_{\substack{\omega\in\D_f \\ \omega<\kappa_0}} A_\omega(\bfalpha) n^{-\omega^*}
\]
with suitable coefficients $A_\omega(\bfalpha)$; we do not enter the properties of such coefficients, which are definitely useful but not easy to describe. With such a notation, the transformation formula reads as follows.

\medskip
{\bf Theorem 3.3.} (Kaczorowski-Perelli \cite{Ka-Pe/2011a},\cite{Ka-Pe/resoII}) {\sl Let $F\in\S$ with $d_F>0$ and $\theta_F=0$, let $f(n,\bfalpha)$ be as in \eqref{3-1} with $\kappa_0>1/d_F$ and $\alpha_0>0$, and let $K\geq0$. There exist an integer $J\geq0$, constants $0=\eta_0<\eta_1<\dots<\eta_J$ and functions $W_0(s),\dots,W_J(s)$, $G_J(s;f)$ holomorphic for $\sigma>-K$, with $W_0(s)$ nonvanishing there, such that for }$\sigma>-K$
\begin{equation}
\label{3-2}
F(s;f) = \sum_{j=0}^J W_j(s) \overline{F}(s^*+\eta_j;f^*) + G_J(s;f).
\end{equation}

\medskip
As in the case $\kappa_0\leq 1/d_F$, there are bounds for the order of growth of $G_J(s;f)$ on vertical lines. The meaning of \eqref{3-2} is that the difference between the left hand side and the sum on the right hand side is holomorphic for $\si>-K$; hence in general this is essentially a relation between two nonlinear twists whose properties are unknown outside their half-planes of absolute convergence. Actually, thanks to the link between $\Re{s}$ and $\Re{s^*}$ one can often get directly from \eqref{3-2} some analytic continuation for $F(s;f)$ (or for $F(s;f^*)$), but apparently not much more. Moreover, the {\it duality operator} $T$, sending $f\mapsto T(f)=f^*$, is self-reciprocal (i.e. $T^2(f)=f$), thus preventing further immediate applications of the transformation formula. Hence another idea is needed in order to exploit the full force of Theorem 3.3.

\medskip
The required idea is actually quite simple, but quite efficient as well. Indeed, given a polynomial $P\in\ZZ[x]$ and $f(n,\bfalpha)$ as in \eqref{3-1}, we consider the {\it shift operator} $S$ defined by
\begin{equation}
\label{3-3}
S(f)(n,\bfalpha) = f(n,\bfalpha) + P(n).
\end{equation}
In view of the periodicity of the complex exponential, $S$ acts trivially on nonlinear twists since $F(s;S(f)) = F(s;f)$, but it turns out that the action of $T$ on $S(f)$ is different from its action on $f$. Hence we consider the noncommuntative group $\frak{G}$ of the operators of type
\[
T S_k T S_{k-1}\dots TS_1,
\] 
where the $S_j$ are shifts as in \eqref{3-3}, thus getting a plenty of new transformation formulae for $F(s;f)$. Of course, before any application of $T$, the choice of $S_j$ must ensure that the new leading exponent is $>1/d_F$. In particular, we may start choosing a nonlinear twist $F(s;f)$ with $0\leq \kappa_0\leq 1/d_F$, whose analytic properties are quite well known, apply a shift $S_1$, which does not change the initial choice but produces a new leading exponent $>1/d_F$, and then apply $T$. From Theorem 3.3 we therefore get an expression of type \eqref{3-2} with $f^*=TS_1(f)$, but now we have information on the analytic properties of $F(s;f)$. Since $\eta_0<\eta_1<\dots<\eta_J$, this provides the meromorphic continuation of $F(s;f^*)$ in a half-plane larger than $\si>1$; hence by an iteration we get the meromorphic continuation of $F(s;f^*)$ over $\CC$. Moreover, we may apply a new shift $S_2$ to $f^*$, then $T$ again, and so on. This procedure has been formalized in \cite{Ka-Pe/resoII}, and leads to the analytic properties of a new class, denoted by $\A(F)$, of nonlinear twists of $F\in\S$ with $\kappa_0>1/d_F$. Moreover, in \cite{Ka-Pe/resoII} one also finds explicit examples involving classical $L$-functions.

\medskip
We finally mention that the results presented in this subsection do not depend on the last two, more arithmetical, axioms of $\S$. Indeed, Theorems 3.1--3.3 hold in the more general analytic setting of the {\it extended Selberg class} $\S^\sharp$, namely the class of functions $F(s)$ satisfying only axioms {\it (i)}, {\it (ii)} and {\it (iii)} above. As a consequence, we have in particular that the results on the degree conjecture presented in the next subsection hold true for $\S^\sharp$ as well. Moreover, contrary to $\S$, the subclass of the degree 0 functions in $\S^\sharp$ contains classical zeta functions and has interesting properties; see Kaczorowski-Molteni-Perelli \cite{K-M-P/2008}.

\bigskip
{\bf 3.3. Converse theorems for the Selberg class.} We already mentioned the converse theorem for degree 0 functions of $\S$, namely

\medskip
{\bf Theorem 3.4.} (Conrey-Ghosh \cite{Co-Gh/1993}) {\sl The only function in $\S$ of degree $0$ is the constant $1$.}

\medskip
Theorem 3.4 is a consequence of the following result in \cite{Co-Gh/1993}: 

\smallskip
\centerline{{\it the Dirichlet series of every $F\in\S$ with $0\leq d_F<1$ is absolutely convergent over $\CC$.}}

\smallskip
\noindent
Moreover, such a result easily implies the degree conjecture for $0<d_F<1$. Actually, the same principle is already present in papers by Richert \cite{Ric/1957} and Bochner \cite{Boc/1958}, although the three proofs are technically different. Hence we have

\medskip
{\bf Theorem 3.5.} (Richert \cite{Ric/1957}, Bochner \cite{Boc/1958}, Conrey-Ghosh \cite{Co-Gh/1993}) {\sl There are no functions $F\in\S$ with $0<d_F<1$.}

\medskip
As remarked in the previous subsection, a simple and direct proof of Theorem 3.5 can be obtained from the properties of the standard twist. Indeed, by Theorem 3.1 with general $\theta_F\in\RR$, the standard twist $F(s,\alpha)$ with $\alpha\in$ Spec$(F)$ has a simple pole on the line $\si=(d_F+1)/(2d_F)$; but $(d_F+1)/(2d_F)>1$ if $0<d_F<1$, giving a contradiction since the Dirichlet series of $F(s,\alpha)$ is absolutely convergent for $\si>1$.

\medskip
The next case are the degree 1 functions of $\S$; the general converse conjecture predicts that these functions coincide with the automorphic $L$-functions of degree 1, which are essentially the Dirichlet $L$-functions with primitive characters. In this case the standard twist coincides with the linear twist $F(s,\alpha)$, which has the important property of being periodic in $\alpha$; the proof of the converse theorem depends on such a periodicity.

\medskip
{\bf Theorem 3.6.} (Kaczorowski-Perelli \cite{Ka-Pe/1999a}) {\sl Suppose that $F\in\S$ has $d_F=1$. Then $F(s) = L(s+i\theta,\chi)$ for some primitive Dirichlet character $\chi$ and $\theta\in\RR$.}

\medskip
Theorem 3.6 therefore shows that the general converse conjecture holds true for degree 1; in particular, if $F(s)$ has a pole at $s=1$ then $F(s)=\zeta(s)$. Note that the periodicity of the normalized coefficients $a(n)n^{-i\theta}$ of the degree 1 functions depends ultimately on the $\alpha$-periodicity of their standard twist. We illustrate the proof of Theorem 3.6 by giving a short proof of a general form of Hamburger's theorem; the proof of Theorem 3.6 is technically definitely more involved, but the basic principle is the same. 

\medskip
First note that the solutions of functional equation \eqref{2-1} have degree and conductor 1, although \eqref{2-1} has a very special shape with respect to the a priori general functional equations satisfied by functions of degree and conductor 1. Of course, strictly speaking \eqref{2-1} is not of Selberg class type since $F(s)$ is not reflected to its conjugate, but this is a minor point. Therefore, we may regard as a general form of Hamburger's theorem the assertion that {\it if a function $F(s)$ from the extended Selberg class $\S^\sharp$ has degree and conductor $1$, then $F(s)=c\zeta(s)$ for some constant $c$}. We recall that the results in the previous subsection hold also for $\S^\sharp$; moreover, since we stated Theorem 3.1 with the additional condition $\theta_F=0$, for simplicity we add this condition here as well. Consider now the standard twist $F(s,\alpha)$, and recall notation and results in Theorem 3.1. In this case we have $s_0=1$ and $n_\alpha=\alpha$, hence choosing an integer $n_0\geq1$ with $a(n_0)\neq0$ we have
\[
\res_{s=1} F(s,n_0) = c_0(F) \overline{a(n_0)} \hskip2cm c_0(F)\neq0.
\]
But from the 1-periodicity in $\alpha$ of $F(s,\alpha)$ we have that $F(s,\alpha)$ has a simple pole at $s=1$ with the same residue $c_0(F)\overline{a(n_0)}$ whenever $\alpha=n$ with $n\in\NN$. Hence, again by Theorem 3.1, we deduce that Spec$(F)=\NN$, and for every $n\in\NN$
\[
c_0(F) \overline{a(n)} = \res_{s=1} F(s,n) =  c_0(F) \overline{a(n_0)}.
\]
Therefore $F(s) = c\zeta(s)$ for some constant $c$. Moreover, if $F(s)$ is an Euler product, then necessarily $c=1$.

\medskip
Now we turn to the degree conjecture in the range $1<d_F<2$.

\medskip
{\bf Theorem 3.7.} (Kaczorowski-Perelli \cite{Ka-Pe/2011a}) {\sl There are no functions $F\in\S$ with $1<d_F<2$.}

\medskip
In this case the properties of the standard twist alone are apparently not sufficient, and indeed the proof involves both the standard twist and Theorem 3.3. The basic strategy is, as for Theorem 3.5, to get a contradiction by inducing a pole of an auxiliary nonlinear twist in its half-plane of absolute convergence; however, the actual realization of such a strategy requires deeper ideas. Again, we illustrate the proof of Theorem 3.7 by means of a special case.

\medskip
Given $F\in\S$ with $1<d_F<2$, if it exists, we consider its standard twist $F(s,\alpha)$ and write $f_0(n,\alpha) = \alpha n^{1/d_F}$, so that $F(s,\alpha) = F(s;f_0)$. Then we apply the shift $S(f_0)(n,\alpha) = f_0(n,\alpha) + n$ and then the operator $T$, thus getting an expression of type \eqref{3-2} with $f=f_0$ and $f^* = TS(f_0)$. Next we compute the real part of the corresponding $s^*$ by means of \eqref{3-1bis}. In this case we have $\kappa_0=1$, hence
\[
\Re{s^*} = \frac{\si+d_F/2-1}{d_F-1}.
\]
Now choose $\alpha\in$ Spec$(F)$, thus $F(s;f_0)$ has a pole at $s=s_0$ with
\[
\Re{s_0} = \frac{d_F+1}{2d_F}
\]
(again we use Theorem 3.1 with general $\theta_F\in\RR$). As a consequence, since $\eta_0<\eta_1<\dots<\eta_N$ and $G_J(s,f_0)$ is holomorphic, $\overline{F}(s^*;f^*)$ has a pole at $s^*=s_0^*$ with
\[
\Re{s_0^*} = \frac{d_F^2-d_F+1}{2d_F(d_F-1)}.
\]
Therefore, $\overline{F}(s^*;f^*)$ has a pole in $\Re{s^*}>1$ if $1<d_F<\frac{1+\sqrt{5}}{2}$, a contradiction showing that there are no functions in $\S$ with $1<d_F<\frac{1+\sqrt{5}}{2}$.

\medskip
In order to prove Theorem 3.7 for the full interval $1<d_F<2$, we need an algorithm to construct iterates of type $TS_kTS_{k-1}\dots TS_1(f_0)$ such that the associated nonlinear  twist has a pole in the absolute convergence half-plane, for every $1<d_F<2$; this is what is actually done in \cite{Ka-Pe/2011a}. Here we only say that such an algorithm is quite complex and its length tends exponentially to $\infty$ as $d_F\to 2$. Moreover, the algorithm is designed for the interval $1<d_F<2$, and at present it is not at all clear how to extend it to higher degrees.

\medskip
We remark in passing that for degrees $<2$ there is a complete theory for the extended Selberg class as well. Indeed, not only the degree conjecture holds for $\S^\sharp$, but also the converse theorems are known for degrees 0 and 1; we refer to \cite{Ka-Pe/1999a} for such results. However, while the degree conjecture is expected to hold in general for $\S^\sharp$, 
as far as we know there is no reasonable guess on the shape of converse theorems for $\S^\sharp$, for integer degrees $\geq2$.

\medskip
The next open case is the degree 2. In this case the general converse conjecture predicts that $\S$ contains, essentially, the $L$-functions associated with holomorphic and nonholomorphic eigenforms. For degree 2 functions, the linear twist is {\it self-dual}, in the sense that application of the $T$-operator to a linear twist produces another linear twist. Moreover, the subcase of the $F\in\S$ with $q_F=1$ offers a concrete possibility of attack thanks to a closer control on the analytic properties of linear twists. In particular, when $F(s)$ has a pole at $s=1$ we have

\medskip
{\bf Theorem 3.8.} (Kaczorowski-Perelli \cite{Ka-Pe/2015}) {\sl Suppose that $F\in\S$ has $d_F=2$, $q_F=1$ and a pole at $s=1$. Then $F(s) = \zeta(s)^2$.}

\medskip
Theorem 3.8 might look quite special, but nevertheless its proof requires new ideas with respect to all previous cases. The main ingredient is the use of analytic properties of the linear twists to get control on the shape of the Euler product of $F(s)$; apparently, such a link is not present in the classical theory of $L$-functions. One may compare Theorem 3.8 with the well known extension by Maass \cite{Maa/1949} of Hamburger's converse theorem, asserting that the solutions $F(s)$ (with the usual analytic properties) of the functional equation satisfied by $\zeta(s)^2$ are all of the form $F(s)=c\zeta(s)^2$ with some constant $c$. Note that Theorem 3.8 and Maass' converse theorem, although both dealing with functions of degree 2 and conductor 1, are of different nature. Indeed, Maass' theorem needs no Euler product but deals with a special functional equation, whose shape allows to get and exploit special properties of certain Mellin transforms. Instead, Theorem 3.8 deals with a general functional equation, thus preventing the possibility of taking advantage from special functions, but the Euler product plays an important role. Actually, the two methods of proof are completely different.

\medskip
We conclude with an outline of work in progress with J.Kaczorowski toward a converse theorem in the general case of functions $F\in\S$ of degree 2 and conductor 1. A first result is in \cite{Ka-Pe/weak}, and may be regarded as a weak form of  converse theorem. Indeed, it is proved that, once normalized by condition $\theta_F=0$, such functions $F(s)$ have real coefficients as in the case of the $L$-functions of holomorphic and nonholomorphic eigenforms of conductor 1; moreover, the first $H$-invariant $H_F(1)$ and the root number also agree with the corresponding invariants of these $L$-functions. The proof is based on a close analysis of the transformation formula \eqref{3-2} with twist functions of type $f(n,\bfalpha) = \alpha_0 n + \alpha_1 \sqrt{n}$. Moreover, there are indications that a sharp converse theorem based on the coincidence of the first few $H$-invariants is within reach in  the degree 2 and conductor 1 case.

\bigskip
{\bf Acknowledgments.} We wish to thank Sandro Bettin, Giuseppe Molteni, Stefano Vigni and the referee for carefully reading the paper and suggesting several improvements. The author is member of the GNAMPA group of INdAM.

\newpage

\ifx\undefined\bysame{poly}.
\newcommand{\bysame}{\leavevmode\hbox to3em{\hrulefill}\ ,}
\fi

\bigskip
\bigskip
Alberto Perelli

Dipartimento di Matematica

Universit\`a di Genova

via Dodecaneso 35

16146 Genova, Italy

e-mail: perelli@dima.unige.it

\end{document}